\documentclass[10pt]{article}
\usepackage{amssymb, amsmath}

\begin{document}
\title{On the number of solutions of a quadratic equation\\ in a normed space}
\author{Victor Alexandrov}
\date{ }
\maketitle
\begin{abstract}
We study an equation $Qu=g$, where $Q$ is a continuous quadratic operator acting from
one normed space to another normed space.
Obviously, if $u$ is a solution of such equation then $-u$ is also a solution.
We find conditions implying that there are no other solutions and
apply them to the study of the Dirichlet boundary 
value problem for the partial differential equation  
$u\Delta u =g$.
\end{abstract}

\bigskip 

\centerline{\S~1. Continuous quadratic operators}

\centerline{and continuous symmetric bilinear operators}

\medskip

Let $L$ and $L'$ be normed spaces over the field 
of real numbers $\mathbb R$.

\textbf{Definition 1:}
A continuous mapping $Q:L\to L'$ is called a 
\textit{continuous quadratic operator} if
$$
Q(u+v)+Q(u-v)= 2Q(u)+2Q(v)
\quad\mbox{for all $u,v\in L$}\quad\eqno(1)
$$
and
$$
Q(ku)= k^2Q(u)
\quad\mbox{for all $u\in L$ and all $k\in\mathbb R$}.\quad\eqno(2)
$$

\textbf{Definition 2:}
Given $B: L\times L\to L'$, a continuous symmetric bilinear 
mapping, the mapping $Q:L\to L'$, defined by the formula
$$
Q(u)=B(u,u),\eqno(3)
$$
is called a \textit{continuous quadratic operator}.

\textbf{Theorem 1:}
\textit{Definitions 1 and 2 are equivalent.}

\textbf{Proof:}
Let $Q$ be as in Definition 1.
By definition, put
$$
B(u,v)=\frac{1}{4}\bigl(Q(u+v)-Q(u-v)\bigr). \eqno(4)
$$

Obviously, $B$ is continuous, $Q$ is generated by $B$ according 
to the formula (3), and $B$ is symmetric. Less trivial is
that $B$ is bilinear. In order to prove the letter, we 
define mappings $F:L\times L\times L\to L'$ 
and $f:\mathbb R\times L\times L\to L'$ by the formulas
$$
F(u,v,w)=4\bigl(B(u+v,w)-B(u,w)-B(v,w)\bigr)
$$
and
$$
f(k)=B(ku,v)-kB(u,v),
$$
where $B$ is defined by (4). Our aim is to prove that 
$$
F(u,v,w)=0 
\quad\mbox{for all $u,v,w\in L$}\quad\eqno(5)
$$
and
$$
f(k,u,v)=0 
\quad\mbox{for all $u,v\in L$ and all $k\in\mathbb R$.}\quad\eqno(6)
$$

As soon as formulas (5) and (6) will be proved,
we can conclude that $B$ is bilinear and, thus, Definition 2 
follows from Definition 1.

In order to prove (5), we use (4) and rewrite $F$ in the form
$$
F(u,v,w)=Q(u+v+w)-Q(u+v-w)-Q(u+w)+Q(u-w)-Q(v+w)+Q(v-w). \eqno(7)
$$
It follows from (1) that
$$
Q(u+v+w)=2Q(u+w)+2Q(v)-Q(u-v+w)\eqno(8)
$$
and
$$
Q(u+v-w)=2Q(u-w)+2Q(v)-Q(u-v-w).\eqno(9)
$$
Substituting (8) and (9) to (7) yields
$$
F(u,v,w)=Q(u-v-w)-Q(u-v+w)+Q(u+w)-Q(u-w)+Q(v-w)-Q(v+w). \eqno(10)
$$
Adding (10) to (7) and taking into account (2), we get
$2F(u,v,w)=0$. Thus (5) is proved.

In order to prove (6), we use (4) and rewrite $f$ in the form
$$
f(k,u,v)=\frac{1}{4}\bigl(Q(ku+v)-Q(ku-v)\bigr)-
\frac{1}{4}\bigl(Q(u+v)-Q(u-v)\bigr).\eqno(11)
$$
From (11), we find $f(0,u,v)=0$ and, using (2), 
$$
f(-1,u,v)=0.\eqno(12)
$$ 
Given any integer $n\neq 0$, we use (5) and (12) 
to obtain
\begin{align*}
f(n,u,v) & =B(nu,v)-nB(u,v)=
B\bigl(\mbox{sgn\,} n(u+u+\dots+u),v\bigr)-nB(u,v)\\
& =\mbox{sgn\,} n \bigl(B(u,v)+B(u,v)+\dots+B(u,v)\bigr)-nB(u,v)=  
nB(u,v) -nB(u,v)= 0.
\end{align*}
Here $\mbox{sgn\,} n$ denotes the sign of $n$, i.e.,
$\mbox{sgn\,} n=+1$ for $n>0$,
$\mbox{sgn\,} n=0$ for $n=0$, and
$\mbox{sgn\,} n=-1$ for $n<0$.

If $k$ is a rational number, $k=m/n$, then
\begin{align*}
f(k,u,v) & =
B\biggl(\frac{m}{n}u,v\biggr)-\frac{m}{n}B(u,v)=
mB\biggl(\frac{1}{n}u,v\biggr)-\frac{m}{n}B(u,v)\\
& =\frac{m}{n}nB\biggl(\frac{1}{n}u,v\biggr)-\frac{m}{n}B(u,v)=
\frac{m}{n}B(u,v)-\frac{m}{n}B(u,v)=0.
\end{align*}
Thus, $f(k,u,v)=0$ for all rational $k$ and all $u,v\in L$.
Since $f$ is continuous, $f(k,u,v)=0$ for all $k\in\mathbb R$ 
and all $u,v\in L$. This completes the proof of (6).
Hence, if $Q$ satisfies Definition 2, it satisfies Definition 1 also.

At last, observe that straightforward computations show that, if $Q$ 
satisfies Definition 1, it satisfies Definition 2 also. 
Theorem 1 is proved.

Theorem 1 shows that there is a one-to-one correspondence
between continuous quadratic operators $Q:L\to L'$
and continuous symmetric bilinear operators $B:L\times L\to L'$.
By $B_Q$ we denote the unique continuous symmetric bilinear operator
that corresponds to a continuous quadratic operator $Q$.
Thus, $Q(u)=B_Q(u,u)$ for all $u\in L$.

Quadratic operators were studied by many authors from
different points of view. For more details the reader is 
referred, for example, to \cite{AC12} and references given therein.
We are going to study a quadratic equation $Qu=g$, where $Q:L\to L'$ is a 
continuous quadratic operator acting from a normed space $L$ to another normed space $L'$.
It follows from (2) that if $u$ is a solution of a quadratic equation then $-u$ 
is also a solution. We are interested in the following question:
given a continuous quadratic operator $Q$, is it true that, for all $g$, 
the equation $Qu=g$ has no more than two solutions?
The answer is given by the following

\textbf{Theorem 2:} 
\textit{Let $L$ and $L'$ be normed spaces over
the field $\mathbb R$ of real numbers, $Q:L\to L'$
be a continuous quadratic operator. 
Then the following conditions are equivalent:}

(i) \textit{for all $g\in L'$, the equation $Qu=g$ has no more then two solutions, 
i.e., $Q(u)=Q(v)$ implies $v=\pm u$};

(ii) \textit{$B_Q$ is nondegenerate, i.e., for every $u\in L$,
$u\neq 0$, the equality $B_Q(u,v)=0$ implies $v=0$.}

\textbf{Proof:}
Suppose that condition (ii) is not satisfied.
This means that, in $L$, there exist $u\neq 0$ and $v\neq 0$ such that $B(u,v)=0$.
Then, according to (4), $Q(u+v)-Q(u-v)=4B_Q(u,v)=0$.
On the other hand, $u+v\neq u-v$ as well as $u+v\neq -(u-v)$. 
Hence, condition (i) is not satisfied. 
In other words, this means that condition (i) implies condition (ii).

Now suppose that condition (i) is not satisfied.
This means that there are two vectors $u,v\in L$ such that $Q(u)=Q(v)$ and neither 
$v=u$ nor $v=-u$. By definition, put $U=\frac12(u+v)$ and $V=\frac12(u-v)$. 
Obviously, $U\neq 0$ and $V\neq 0$. 
Moreover, according to (4), $4B_Q(U,V)=Q(U+V)-Q(U-V)=Q(u)-Q(v)=0$.
Hence, condition (ii) is not satisfied. 
In other words, condition (ii) implies condition (i).
Theorem 2 is proved.

\bigskip 

\centerline{\S~2. On the number of solutions to a boundary 
value problem}

\centerline{for the partial differential equation 
$u\Delta u =g$}

\medskip

In order to demonstrate the advantages of the
replacement of the quadratic equation $Qu=g$
by the linear equation $B_Q(u,v)=0$ ($g$ and $v$ are prescribed functions)
described in the previous section, we consider the following Dirichlet problem:
\begin{align*}
&u(x,y)\biggl(
\frac{\partial^2 u}{\partial x^2}(x,y)+\frac{\partial^2 u}{\partial y^2}(x,y)
\biggr)=g(x,y) \quad\textrm{for}\ (x,y)\in D, \tag{13}\\
&u(x,y)=0 \quad\textrm{for}\ (x,y)\in \partial D, \tag{14}
\end{align*}
where 
$$D=\biggl\{(x,y)\in\mathbb R^2\biggl\vert\biggr. -\frac{\pi}{2}<x<\frac{\pi}{2}, 
-\frac{\pi}{2}<y<\frac{\pi}{2}\biggr\}
$$
is a square and $\partial D$ is its boundary.

Obviously, if $u$ is a solution to (13)--(14) then $-u$ is also a solution.  
We are interested if, for every $f$, the problem (13)--(14) has no more than two solutions.
Equivalently, our problem can be reformulated as follows:
is it true that if $u$ and $v$ are two solutions to (13)--(14) then $v=\pm u$?

In this section we prove that the answer is negative.

The nonlinear PDE operator in the left-hand side of (13) is a simplification
of the operator in the equation 
$u\Delta u+C_1|\nabla u|^2=C_2$, where $C_1, C_2\in\mathbb R$,
for which a strong minimum principle is proved among other results in \cite{Ma02}.

Let $L$ be the linear space of functions $u:D\to\mathbb R$, each of which is of 
the class $C^{\infty}$ in some (specific for every $u$) open set containing 
the closure of $D$ and vanishes on $\partial D$, endowed with the norm
$\| u\|_L=\sup_{(x,y)\in D}|u(x,y)|+\sup_{(x,y)\in D}|\nabla u(x,y)|$,
where
$$
\nabla u(x,y)=\biggl(
\frac{\partial u}{\partial x}(x,y), \frac{\partial u}{\partial y}(x,y)
\biggr)
$$
is the gradient of $u$. Let $L'$ be the linear space of functions which are continuous
in the closure of $D$ and are endowed with the norm
$\| u\|_{L'}=\sup_{(x,y)\in D}|u(x,y)|$.
At last, let the operator $Q:L\to L'$ be defined by the formula
$Qu=u\Delta u$ or, equivalently,
$$
(Qu)(x,y)=
u(x,y)\biggl(
\frac{\partial^2 u}{\partial x^2}(x,y)+ 
\frac{\partial^2 u}{\partial y^2}(x,y)
\biggr).
$$

Obviously, $Q$ is a continuous quadratic operator from $L$ to $L'$ for which 
$B_Q(u,v)=u\Delta v+ v\Delta u$
and the equation $Qu=g$ coincides with (13).

\textbf{Theorem 3:} 
\textit{Let $D$, $L$, $Q$, and $B_Q$ be as defined above in this section.
Then there exist functions $u, v\in L$ that satisfy the equation $B_Q(u,v)=0$ and
do not vanish identically in $D$.
}

\textbf{Proof:}
Obviously, the equation $u\Delta v+ v\Delta u=0$ will be satisfied if there exists
a function $p$ such that
\begin{align*}
&\Delta u(x,y)+p(x,y)u(x,y)=0, \tag{15}\\
&\Delta v(x,y)-p(x,y)v(x,y)=0 \tag{16}
\end{align*}
for $(x,y)\in D$. 

We look for solutions to (15)--(16) of a special form, namely,
$u(x,y)=U(x)\widetilde{U}(y)$ and
$v(x,y)=V(x)\widetilde{V}(y)$.
In order to satisfy the boundary condition (14), we assume that
$$
U\biggl(\pm\frac{\pi}{2}\biggr)=\widetilde{U}\biggl(\pm\frac{\pi}{2}\biggr)=
V\biggl(\pm\frac{\pi}{2}\biggr)=\widetilde{V}\biggl(\pm\frac{\pi}{2}\biggr)=0.\eqno(17)
$$
Moreover, we assume that $p$ is a function of a single variable $x$,
i.\,e., that $p$ is independent of $y$. Under these assumptions, (15)--(16) yield
\begin{align*}
&U''(x)\widetilde{U}(y)+U(x)\widetilde{U}''(y)+p(x)U(x)\widetilde{U}(y)=0, \\
&V''(x)\widetilde{V}(y)+V(x)\widetilde{V}''(y)-p(x)V(x)\widetilde{V}(y)=0,
\end{align*}
or, after separation of variables, 
\begin{align*}
&\frac{U''(x)}{U(x)}+p(x)=-\frac{\widetilde{U}''(y)}{\widetilde{U}(y)}=\lambda, \tag{18}\\
&\frac{V''(x)}{V(x)}-p(x)=-\frac{\widetilde{V}''(y)}{\widetilde{V}(y)}=\mu, \tag{19}
\end{align*}
where $\lambda$ and $\mu$ are some real constants. 
For our purpose it is sufficient to find at least one non-zero solution 
to (18) and (19) satisfying the boundary condition (17). 
So, we put $\lambda=\mu=1$ for simplicity.
Then we easily check that $\widetilde{U}(y)=\widetilde{V}(y)=\cos y$ satisfy
both the equations (18)--(19) and boundary conditions (17). 

The problem of finding the functions $U(x)$, $V(x)$, and $p(x)$ satisfying
the equations (18)--(19) and boundary conditions (17) is more complicated.
We put by definition $a=-1$ and $p(x)=-2q\cos 2x$, where $q$ is a real number
that will be specified later. Then the equation (18) implies that
the function $U(x)$ satisfies the equation
$$
z''(x)+(a-2q\cos 2x)z(x)=0.\eqno(20)
$$
Similarly, the equation (19) implies that
the function $V(x)$ satisfies the equation
$$
z''(x)+(a+2q\cos 2x)z(x)=0,\eqno(21)
$$
which can be obtained from (20) by substituting $-q$ instead of $q$.

The equation (20) is known as the Mathieu equation.
Properties of its solutions are studied in hundreds of articles and books
among which we mention only three classical treatises \cite{St32}, \cite{Mc47}, and \cite{An52}
and the chapters \cite{Bl64}, \cite{Wo10} in well-known handbooks.

Recall from \cite{St32}--\cite{Wo10} 
that one of periodic solutions of the Mathieu equation (20) 
is usually denoted by $\textrm{se}_2(x,q)$ and may be 
defined by the following expansion
$$
\textrm{se}_2(x,q)=\sum\limits_{m=1}^{\infty}B_{2m}^2 \sin(2mx),\eqno(22)
$$
where the coefficients $B_{2m}^2$ satisfy the recurrence relations
$qB_4^2=(b_2-4)B_2^2$ and $qB_{2m+2}^2=(b_2-4m^2)B_{2m}^2-qB_{2m-2}^2$ for $m\geqslant 2$,
as well as the following normalization condition 
$\sum_{m=1}^{\infty} (B_{2m}^2)^2=1$.

Recall also that, for a given value of $q$, the function $\textrm{se}_2(x,q)$ 
is a solution of the Mathieu equation (20) not for arbitrary values of the parameter $a$;
the parameter $a$, referred to as an eigenvalue, depends on $q$ and, in fact, 
is a function of $q$. Usually,  this  function is denoted by $a= b_2(q)$. 
Though this function is quiet complicated, it is known that
$b_2(-q)=b_2(q)$ and there is $q_*\approx 8$ such that $b_2(q_*)=-1$. 

From these properties of the function $\textrm{se}_2(x,q)$ we immediately conclude that
the function $U(x)=\textrm{se}_2(x,q_*)$ satisfies the equation (20) and the boundary 
conditions (17). 
We conclude also that the function $V(x)=\textrm{se}_2(x,-q_*)$ satisfies the equation (21) 
and the boundary conditions (17). 

Hence, the functions $u(x,y)=\textrm{se}_2(x,q_*)\cos y$ and 
$v(x,y)=\textrm{se}_2(x,-q_*)\cos y$ satisfy the equation $u\Delta v+ v\Delta u=0$ in $D$
and vanish on $\partial D$. This completes the proof of Theorem 3.

\textbf{Theorem 4:} 
\textit{Let $D$, $L$, and $Q$ be as defined at the beginning of this section.
Then there exist functions $u, v\in L$ that satisfy the equation $Q(u)=Q(v)$ and
neither $u=v$, nor $u=-v$.
}

\textbf{Proof} follows immediately from Theorems 2 and 3.

\bigskip

\noindent{Victor Alexandrov}

\noindent\textit{Sobolev Institute of Mathematics}

\noindent\textit{Koptyug ave., 4}

\noindent\textit{Novosibirsk, 630090, Russia}

and

\noindent\textit{Department of Physics}

\noindent\textit{Novosibirsk State University}

\noindent\textit{Pirogov str., 2}

\noindent\textit{Novosibirsk, 630090, Russia}

\noindent\textit{e-mail: alex@math.nsc.ru}

\bigskip

\noindent{Submitted: June 8, 2015}

\end{document}